\documentclass[11pt]{article}
\usepackage{amsmath,amsthm}

\def\ftoday{le \space\number\day \space\ifcase\month\or
  janvier\or f\'evrier\or mars\or avril\or mai\or juin\or
  juillet\or ao\^ut\or septembre\or octobre\or novembre\or d\'ecembre\fi
  \space\number\year}

\def\real{I\kern-0.20em R}

\def\integer{I\kern-0.20em N}
\def\relative{{\rm \rlap Z\kern 2.2pt Z}}
\def\cc{\kern-.25em{\c c}}

\def\bc{\begin{center}}
\def\ec{\end{center}}
\def\=def{\stackrel{{\rm def}}{=}}

\newcounter{indconst}
\newcounter{auxconst}

\def\bit{\begin{itemize}}
\def\eit{\end{itemize}}
\def\ben{\begin{enumerate}}
\def\een{\end{enumerate}}
\def\bde{\begin{description}}
\def\ede{\end{description}}
 
\def\beq{\begin{equation}}
\def\eeq{\end{equation}}
\def\bfi{\begin{figure}[hbt] \begin{center}}
\def\efi{\end{center} \end{figure}}

\def\bce{\begin{center}}
\def\ece{\end{center}}

\newcommand{\Proof}{\medskip\noindent{\underline{Proof:}}\quad\medskip}        

\newtheorem {theos} {Theorem}[section]

\newtheorem {coros} {Corollary} [section]

\newtheorem {lemms} {Lemma}[section]

\newtheorem {props} {Proposition}[section]

\newtheorem {defis} {Definition}[section]

\newtheorem {rems}{Remark} [section]
\newtheorem {ex} {Exemple}[section]

\input amssym.def
\input amssym
\begin{document}

\title{Small divisors and large multipliers}
\author{Boele Braaksma\thanks{University of Groningen, Department of Mathematics
P.O. Box 800, 9700 AV Groningen, The Netherlands. E-mail : {\tt braaksma@rug.nl}} \hspace{1mm} and Laurent Stolovitch\thanks{CNRS UMR 5580, 
Universit\'e Paul Sabatier, MIG, Laboratoire de Math\'ematiques Emile Picard,
31062 Toulouse cedex 9, France. E-mail : {\tt stolo@picard.ups-tlse.fr}} }
\date{\ftoday}
\maketitle
\def\abstractname{R\'esum\'e}
\begin{abstract}
Nous \'etudions des germes de champs de vecteurs holomorphes singuliers \`a l'origine de $\Bbb C^n$ dont la partie lin\'eaire est $1$-r\'esonante et qui admettent une forme normale polynomiale. En g\'en\'eral, bien que le diff\'eomorphisme formel normalisant soit divergent \`a l'origine, il existe n\'eanmoins des diff\'eomorphismes holomorphes dans des "domaines sectoriels" qui les transforment en leur forme normale. Dans cet article, nous \'etudions la relation qui existe entre le ph\'enom\`ene de petits diviseurs et le caract\`ere Gevrey de ces diff\'eomorphismes sectoriels normalisants. Nous montrons que l'ordre Gevrey de ce dernier est reli\'e au type diophantien des petits diviseurs.
\end{abstract}
\def\abstractname{Abstract}
\begin{abstract}
We study germs of singular holomorphic vector fields at the origin of $\Bbb C^n$ of which the linear part is $1$-resonant and which have a polynomial normal form. The formal normalizing diffeomorphism is usually divergent at the origin but there exists holomorphic diffeomorphisms in some "sectorial domains" which transform these vector fields into their normal form. In this article, we study the interplay between the small divisors phenomenon and the Gevrey character of the sectorial normalizing diffeomorphisms. We show that the Gevrey ordrer of the latter is linked to the diophantine type of the small divisors.
\end{abstract}
\tableofcontents
\section{Introduction}

In this article, we are concerned with the study of some germs of holomorphic vector fields in a neighborhood of a fixed point. More precisely, we shall consider holomorphic non-linear perturbations
$$
X=\sum_{i=1}^n{\left(\lambda_ix_i+f_i(x)\right)\frac{\partial}{\partial x_i}}
$$
of the diagonal linear vector field $s=\sum_{i=1}^n\lambda_ix_i\frac{\partial}{\partial x_i}$, where $n\geq 2$. Hence,
the functions $f_i$'s vanish as well as their first derivatives at the origin. Two such germs of vector fields $X_1, X_2$ are holomorphically conjugate (resp. equivalent) if there exists a germ of biholomorphism $\Phi$ of $(\Bbb C^n,0)$ which conjugates them (resp. up to the muliplication by an holomorphic unit): $\Phi_*X_1(y):=D\Phi(\Phi^{-1}(y))X_1(\Phi^{-1}(y))=X_2(y)$. It is well known (see \cite{Arn2} for instance) that such a vector field is formally conjugate (i.e. by  means of a formal diffeomorphism named {\it normalizing diffeomorphism}) to a {\it normal form} $\widehat X_{norm}$, that is a formal vector field which commutes with the linear part $s$. In coordinates, we have 
$$
\widehat X_{norm}=\sum_{i=1}^n{\left(\lambda_iy_i +  \sum_{
{\scriptstyle (Q,\lambda)=\lambda_i  \atop
\scriptstyle Q \in \Bbb N^n_2 }}
{a_Q^{i}y^Q}\right)\frac{\partial}{\partial y_i}} 
$$
the sum being taken over the multiindices $Q=(q_1,\ldots,q_n)\in \Bbb N^n$ such that  $|Q|:=q_1+\cdots+q_n\geq 2$ (we shall write $Q\in\Bbb N_2^n$) and which satisfy a {\it resonance relation} $(Q,\lambda):=q_1\lambda_1+\cdots+q_n\lambda_n=\lambda_i$. The $a_Q^i$'s are complex numbers.

If there exists a monomial $x^r$ which is a first integral of the linear part $s$ (i.e. $s(x^r)=0$) but which is not a first integral of a normal form then the formal normalizing diffeomorphism is generally a (vector of) divergent power series (see \cite{bruno, Mart1}). We are interested in this situation. 
More precisely, Ichikawa \cite{Ichikawa} has shown that, if $s$ is {\it $1$-resonant} (i.e. the formal non-linear centralizer of $s$ is generated by the sole relation $(r,\lambda)=0$ for a nonzero $r\in \Bbb N^n$), then $X$ has only a finite number of formal invariants if and only if $y^r$ is not a first integral of a normal form (and thus of any normal form); this is Ichikawa's condition I. This means that $X$ has a polynomial normal form.

The analytic classification of these objects in dimension two is due to J. Martinet and J.-P. Ramis in two seminal articles \cite{Ram-Mart1, Ram-Mart2}. They showed that the divergent normalizing diffeomorphism is in fact {\it summable} in some sectorial domain. This means that there exist germs of a holomorphic diffeomorphism in some {\it large sectorial domain} (with vertex at the origin) having the formal diffeomorphism as asymptotic expansion in the domain and conjugating the vector field to its polynomial normal form. The counterpart for one dimensional diffeomophisms is due to Ecalle-Voronin-Malgrange \cite{Voronin, Malg-diffeo}. This study was continued in a more general setting (in particular in higher dimension) by J. Ecalle by a completely different manner in the article \cite{Ecalle-mart}. The second author has given a unified approach of the two articles of Martinet and Ramis while treating the $n$-dimensional case \cite{Stolo-classif}. In this situation, he also proved that the formal normalizing diffeomorphism could be realized as the asymptotic expansion of some genuine holomorphic diffeomorphism in some sectorial domain which normalizes the vector field (in the case of one zero eigenvalue, in dimension $3$ and without small divisors, the result is due to
 Martins \cite{martins-ams}). Due to the presence of small divisors, the summability property of the formal power series does not hold as we already noticed in  \cite{Stolo-classif} (there is no small divisor in the two-dimensional problem). But, we only had a qualitative approach of this phenomenon.
In his article \cite{Ecalle-mart}, J. Ecalle gave some statements (see propositions 9.1-9.4, p.136-137) in a more general setting with {\it "preuves succintes"}. One of them is (we refer to his article for the definitions, $\tilde g(z)$ is a normalizing series):
{\it \begin{quote}
"Th\'eor\`eme 9.1 - Les s\'eries formelles $\tilde g(z)$ associ\'es \`a l'objet local $X$ ou $f$ sont toujours de classe Gevrey $1+\delta^{int}+0$ mais pas inf\'erieur en g\'en\'eral."
\end{quote}}

In this article, we shall {\bf quantify the interplay} between the {\bf Gevrey property} (this is due to the "grands multiplicateurs" that Poincar\'e mentions in \cite{Poincare-methode2}[p.392]) of the normalizing tranformations and the {\bf small divisors phenomenon} (compare with theorem \ref{main2}).\\

{\bf Acknowledgment.}
The authors would like to thank Bernard Malgrange for his remarks and comments that permitted us to improve the redaction of this article. We would like to thank Jean-Pierre Ramis for his encouragement and also for giving us the exact reference in the work of Henri Poincar\'e. We also would like to thank Yann Bugeaud for enlightening discussions and references about arithmetics. The first author would like to thank the department of mathematics of the Universit\'e Paul Sabatier
for several invitations. 

\section{Normal form of $1$-resonant vector fields}

Let $s=\sum_{i=1}^n \lambda_i x_i\partial/\partial x_i$ be a linear diagonal vector field.
\begin{defis}
We shall say that $s$ is $1$-resonant if there is a monomial $x^r$ where 
$r\in \Bbb N_1^n$ such that
if $i\in \{1,\dots ,n\}$, then for all monomials $x^Q$ such that $[s, x^Q\frac{\partial}{\partial x_i}]=0$, we have $x^Q=(x^r)^lx_i$ for some nonnegative integer $l$.
\end{defis}
The resonance monomial $x^r$ generates the ring ${\widehat  {\cal O}_n}^s=\{f\in \Bbb C[[x_1,\ldots, x_n]]\;|\; {\cal L}_s(f)=0 \}$ of first integrals. 
We shall assume that $s$ is diophantine in the sense that it satisfies the
 Bruno small divisors condition 
$$
(\omega)\quad-\sum_{k\geq 0}{\frac{\ln(\omega_{k+1})}{2^k}}<+\infty,
$$
where
$$
\omega_k=\inf\left\{|(Q,\lambda)-\lambda_i|\;|\;|(Q,\lambda)-\lambda_i|\neq 0,\; i=1,\ldots,n,\; Q\in \Bbb N^n,\;2\leq |Q|\leq 2^k\right\}.
$$
Without any loss of generality, we assume that $r_i\neq 0$ if $1\leq i\leq p$ and $r_{p+1}=\cdots=r_n=0$ where $r=(r_1,\ldots, r_n)$. We shall set $p=n$ if none of the $r_i$'s vanish.
Let us define the map $\pi : \Bbb C^n\rightarrow \Bbb C$ defined by $\pi(x)=x^r$. Let $D_{\pi}$ (resp. $\widehat D_{\pi}$) be the group of germs of biholomorphisms (resp. formal diffeomorphisms) at the origin of $\Bbb C^n$, fixing the origin and leaving the map $\pi$ invariant (i.e. $\pi\circ\Phi=\pi$).

Let $k\geq 1$ be an integer and let us define the space ${\cal E}_k$ (resp. $\widehat {\cal E}_k$) of germs of 1-resonant holomorphic (resp. 1-resonant formal) vector fields in a neighboorhood of the origin $0\in \Bbb C^n$ of the form :
$$
X=\sum_{i=1}^n{x_i(\lambda_i+P_{i,k}(x^r))\frac{\partial}{\partial x_i}}+\sum_{i=1}^p{x_i\left(x^r\right)^kf_i(x)\frac{\partial}{\partial x_i}}+\sum_{i=p+1}^n{\left(x^r\right)^{k+1}g_i(x)\frac{\partial}{\partial x_i}}
$$
where $x^r$ is the resonance monomial and $P_{i,k}$'s are polynomials in the variable $u$, vanishing at zero and of degree at most $k$ such that
$\sum_{i=1}^p{r_iP_{i,k}(u)}=\beta u^{k}$ with $\beta\in \Bbb C^*$. The $f_i$'s and $g_i$'s are germs of holomorphic functions (resp. formal power series) in a neighbourhood of $0$ such that $\sum_{i=1}^p{r_if_i(x)}=0$.
We shall say that elements of ${\cal E}_k$ (resp. $\widehat {\cal E}_k$) are {\bf well prepared} vector fields.
We recall that a vector field $X$ is holomorphically equivalent to $Y$ if it is holomorphically conjugate to $Y$ up to multiplication by a unit of ${\cal O}_n$.
\begin{props}\cite{Stolo-classif}[proposition 3.2.1]
Any germ of a $1$-resonant vector field satisfying the Ichikawa transversality condition $(I)$ and $(\omega)$ is holomorphically equivalent to a well prepared germ. This means that there exists an integer $k$ such that $X$ is equivalent to an element of ${\cal E}_k$.
\end{props}

\begin{props}\cite{Stolo-classif}[proposition 3.2.2]
Let $\widehat X\in \widehat {\cal E}_k$ be a well prepared formal vector field of the form
$$
\widehat X=\sum_{i=1}^n{x_i(\lambda_i+P_{i,k}(x^r))\frac{\partial}{\partial x_i}}+\sum_{i=1}^p{x_i\left(x^r\right)^kf_i(x)\frac{\partial}{\partial x_i}}
+\sum_{i=p+1}^n{\left(x^r\right)^{k+1}g_i(x)\frac{\partial}{\partial x_i}}.
$$
Then there exists a unique formal diffeomorphism $\hat \phi\in \widehat D_{\pi}$ 
tangent to the identity at zero such that
$$
\hat \phi_*\widehat X=\sum_{i=1}^n{y_i(\lambda_i+P_{i,k}(y^r))\frac{\partial}{\partial y_i}}.
$$
\end{props}

Let $\alpha\in \Bbb C^n$ such that $(r,\alpha)\neq 0$. Let ${\cal E}_{k,\lambda,\alpha}\subset {\cal E}_k$ be the set of germs of well prepared holomorphic vector fields at the origin of the form:
\beq\label{champs-2}
\sum_{i=1}^n{\left(x_i(\lambda_i+\alpha_{i}(x^r)^k)+(x^r)^kf_i(x)\right)\frac{\partial}{\partial x_i}},
\eeq
with $\sum_{i=1}^p{r_i\alpha_i}=:\beta \neq 0$ and $f_i(x)=x_i\tilde f_i(x)$, $1\leq i\leq p$, where the $\tilde f_i$'s are germs of holomorphic functions in a neighbourhood of $0\in \Bbb C^n$. Moreover, they satisfy $\sum_{i=1}^p{r_i\tilde f_i(x)}=0$. 
Let 
$$
X_{k,\lambda,\alpha}=\sum_{i=1}^n{x_i\left(\lambda_i+\alpha_{i}(x^r)^k\right)\frac{\partial}{\partial x_i}}
$$
be the normal form of such a vector field.

Let us make the following assumptions:
\begin{itemize}
\item $(H'_1)$ all the eigenvalues $\lambda_i$ have a nonnegative imaginary
 part and if zero is not an eigenvalue then there are at least two real
 eigenvalues. 
\item $(H'_2)$ there exists $\lambda_{i_0}\in \Bbb R^*$ such that

$$
\min_{i\neq i_0}{\mbox{Re }\left(\frac{\alpha_i}{\beta}-\frac{\lambda_i}{\lambda_{i_0}}\frac{\alpha_{i_0}}{\beta}\right)}>0,
$$
where the minimum is taken over all indices $i\neq i_0$ such that $\lambda_i
\in \Bbb R.$   
\item $(H'_3)$ if $\lambda_i$ is not real then $f_i=x_i\tilde f_i$ where
 $\tilde f_i$ is a germ of a holomorphic function at the origin.
\end{itemize}

We shall call {\it sectorial domain} a domain of $\Bbb C^n$ of the form:
$$
DS_j(\rho,R)=\left \{y\in \Bbb C^{n}\;\left|\right.\; \left|\arg y^r -\frac{1}{k}\pi(j+\frac{1}{2})\right|<\frac{\pi}{k} -\epsilon, 0<|y^r|<\rho, |y_i|<R\mbox{ for }i=1,\ldots ,n\right \}
$$
where $\rho,R>0, 0<\epsilon <\pi /k$ and $0\leq j\leq 2k-1$ is an integer. 

Stolovitch has shown the following result
\begin{theos}[Sectorial normalization]\cite{Stolo-classif}[Th\'eor\`eme 3.3.1]\label{main-stolo-classif}
Let $\epsilon<\pi/k$ be a positive number and let $X$ belong to ${\cal E}_{k,\lambda,\alpha}$. 
Under assumptions $(H'_1)$, $(H'_2)$ and $(H'_3)$, for any even integer $0\leq j\leq 2k-1$ there exists a local change of coordinates $x_i=y_i+\phi_i^j(y)$, $i=1,\ldots,n$, tangent at the identity, holomorphic in the sectorial domain
$DS_j(\rho,R)$ with $\rho, R$ sufficiently small, in which the vector field $(\ref{champs-2})$ can be written as
\begin{equation}\label{normal-form}
X=\sum_{i=1}^n{y_i\left(\lambda_i+\alpha_{i}(y^r)^k\right)\frac{\partial}{\partial y_i}}.
\end{equation}
This change of coordinates preserves the function $x^r$. 
Each function $\phi_i^j$ admits the formal power series $\hat\phi_i$ as asymptotic expansion in $y^r$ in the sense of G\'erard-Sibuya in the domain 
$DS_j(\rho,R)$. Here, $x_i=y_i+\hat\phi_i(y)$, $i=1,\ldots,n$, is the 
unique formal coordinate system in which the vector field $X$ is in its normal form $(\ref{normal-form})$. Moreover, if all the eigenvalues are real then the 
result holds also for $j$ odd.
\end{theos}

We refer to definition \ref{gerard-sibuya} in the next section for the notion of asymptotic expansion in the sense of G\'erard-Sibuya.
The proof of this theorem reduces to the proof of the sectorial linearization of the non-linear system with an irregular singularity at the origin
\beq\label{irreg}
\beta z^{k+1}\frac{dx_i}{dz}=x_i(\lambda_i+\alpha_iz^k) + z^kf_i(x)\;\;\;\;i=1,\ldots,n.
\eeq
By sectorial linearization, we mean that there is a change of coordinates $x_i=y_i+g_i(z,y)$, $i=1,\ldots, n$, holomorphic in $S\times P$ , where $S$ is a sector in $\Bbb C$ with vertex at $0$ (variable $z$) and $P$ a polydisc centered at $0 \in \Bbb C^n$ (variables $y$) in which the system can be written as:
\beq\label{irreg-lin}
\beta z^{k+1}\frac{dy_i}{dz}=y_i(\lambda_i+\alpha_iz^k)\;\;\;\;i=1,\ldots,n.
\eeq
In \cite{Stolo-classif}) it is shown that the function $\phi_i(y)$ is nothing but $g_i(y^r,y)$. 
Moreover, the $g_i$'s have an expansion at the origin of the form
\beq\label{exp g}
g_i(z,y)=\sum_{Q\in \Bbb N_2^n} g_{i,Q}(z)y^Q
\eeq
where the $g_{i,Q}$'s are holomorphic functions in $S$. 
\section{Gevrey functions and summability}
Here we recall some definitions of Gevrey asymptotics and summability which will 
be used further on (for more details see \cite{malgrange-cours, balser-book, ramis-panora, ramis-stolo-cours}).
\begin{defis}
A holomorphic function $f$ in an open bounded sector $S$ with vertex 0 in 
$\Bbb C$ is said to to admit an asymptotic expansion of Gevrey order $s>0$ (resp. in the sens of Poincar\'e)
if 
there exists a formal power series 
$\hat{f}=\sum _{j=0}^\infty f_j z^j$  such that for every compact subsector
 $S'$ of $S\cup \{0\}$ 
 there exist positive constants $A$ and $C$ such that for all $z\in S'$ and 
$N\in \Bbb N$
$$|f(z)-\sum_{j=0}^{N-1}f_j z^j|\leq CA^N \Gamma(1+Ns)|z|^N,$$
where $\Gamma(x)$ is the Gamma-function (resp. $CA^N \Gamma(1+Ns)$ is to be replaced by a unprecised constant $M_N$).
Such a function $f$ will be called an $s$-Gevrey function on $S$ or shortly
$f$ is $s$-Gevrey on $S$.
\end{defis}
Equivalently: a holomorphic function $f$ on $S$ is $s$-Gevrey on $S$ if
 all derivatives of $f$ are continuous at 0 
and if $S'$ is as above then there exist  positive constants
$A$ and $C$ such that for all $N\in \Bbb N$ and all $z\in S'$:
$$\frac{1}{N!}|\frac{\partial^N f(z)}{\partial z^N}|\leq CA^N\Gamma(1+Ns).
$$
\begin{defis}\cite{Ram-Mart1,Ger-Sib}\label{gerard-sibuya}
Let $\hat f=\sum_{Q\in \Bbb N^n}{\hat f_Q(z)x^Q}\in \Bbb C[[z,x_1,\ldots,x_n]]$ be a formal power series.
We shall say that an analytic function $f$ on $S\setminus \{0\}\times \Delta$ 
( $S\subset \Bbb C$ is an open sector and $\Delta\subset \Bbb C^n$ an open polydisc centered at $0\in \Bbb C^n$), $f(z,x)=\sum_{Q\in \Bbb N^n}{f_Q(z)x^Q}$ admits $\hat f$ as {\bf asymptotic expansion in the sense of G\'erard-Sibuya} in $S\setminus \{0\}\times \Delta$, if each function $f_Q(z)$ admits $\hat f_Q(z)$ as an asymptotic expansion in the sense of Poincar\'e, in the sector $S$ and
for every compact subsector $S'$ of $S\cup \{0\}$, every compact subset $\Delta '$ of $\Delta$ and every $N\in \Bbb N_1$ there exists a constant $K$ such that
$$ |f(z,x)-\sum_{|Q|<N} f_Q(z)x^Q|\leq K|x|^N \mbox{ for all } (z,x)\in S'\times\Delta'.$$

\end{defis}
\begin{defis}
If $k>0$ and $f$ is $1/k$-Gevrey in a sector $S$ with opening $>\pi/k$ then
 $\hat{f}$ is 
 $k$-summable in the direction of the bisector of $S$ 
 and its $k$-sum on $S$
is $f$. In this case $f$ is uniquely determined by $\hat{f}$ and we say that
 $f$ is a $k$-sum on $S$. 
\end{defis}
The notion of summability is due Borel and generalized by Ramis (cf. \cite{ramis-ksum}).
\par Suppose that $f$ is a holomorphic function on $S\times P_n(0,r)$ 
(where $S$ is as above and $P_n(0,r)$ is the open polydisc in ${\Bbb C}^n$ 
with center 0 and radius $r$). Then {\em $f(z,x)$ is said to be $s$-Gevrey in 
$z$ on $S$ uniformly in $x$ on $P_n(0,r')$} for some $r'\in (0,r)$ if 
there exists a formal power series 
$\hat{f}(z,x)=\sum _{j=0}^\infty f_j(x) z^j$ where the coefficients $f_j(x)$
 are holomorphic on 
$P_n(0,r')$ such that for every compact subsector
 $S'$ of $S\cup \{0\}$ 
 there exist positive constants $A$ and $C$ such that for all $z\in S'$, all 
$x\in P_n(0,r')$ and all $N\in \Bbb N$:
\beq\label{unifG}
|f(z,x)-\sum_{j=0}^{N-1}f_j(x) z^j|\leq CA^N \Gamma(1+Ns)|z|^N.
\eeq
The latter condition is equivalent to the existence of $A,C,r'$ as above such
 that for all $z\in S'$, all 
$x\in P_n(0,r')$ and all $N\in \Bbb N$:
$$\frac{1}{N!}|\frac{\partial^N f(z,x)}{\partial z^N}|\leq CA^N\Gamma(1+Ns).
$$
If the opening of $S$ is $>\pi s$ then $\hat{f}(z,x)$ is $1/s$-summable in 
$z$ on $S$ uniformly on $P_n(0,r')$ (cf. \cite{sibuya04}).
\begin{rems}\label{rems1}
If $f(z,x)$ is $s$-Gevrey in $z$ on $S$ uniformly in $x$ on $P_n(0,r')$
then $f(z,x)$ is also $s$-Gevrey in $(z,x)$ on $S\times P_n(0,r'')$ for some
$r''\in (0,r')$ in the sense that for every $S'$ as above there exist positive 
constants $A$ and $C$ such that for all $N\in \Bbb N$ and all $(z,x)\in
S'\times P_n(0,r'')$:
\beq\label{Gevreyzx}
 |f(z,x)-f^{(N-1)}(z,x)|\leq C A^N \Gamma (1+Ns)||(z,x)||^N. 
\eeq
Here $f^{(N-1)}(z,x)$ denotes the $(N-1)$-jet of $f$ in $0\in {\Bbb C}^{n+1}$
and $||.||$ denotes an arbitrary norm on ${\Bbb C}^{n+1}.$
\end{rems}
\Proof We have $f_j(x)=\frac{1}{j!}\frac{\partial^j}{\partial z^j}f(0,x)$ and 
therefore $|f_j(x)|\leq CA^j\Gamma(1+js)$ for $x\in P_n(0,r')$. From this and
 Cauchy's formula it follows that for sufficiently small $\delta>0$ and all
$Q\in {\Bbb N}^n$:
$$\frac{1}{Q!}|\frac{\partial^Q}{\partial x^Q}f_j(0)|\leq CA^j\Gamma(1+js)
\delta^{-|Q|}.$$
Let
$$\tilde{R}_N(z,x)=\sum_{j=0}^{N-1}\sum_{|Q|\geq N-j}\frac{1}{Q!}
\frac{\partial^Q}{\partial x^Q}f_j(0)z^j x^{Q}.$$
Then it follows from (\ref{unifG}) that
\beq\label{tale}
|f(z,x)-f^{(N-1)}(z,x)-\tilde{R}_N(z,x)|\leq CA^N \Gamma(1+Ns)|z|^N
\eeq 
and
$$|\tilde{R}_N(z,x)|\leq C \sum_{j=0}^{N-1}\sum_{|Q|\geq N-j}
A^j\Gamma(1+js)\delta^{-|Q|}|z^j x^Q|. $$
Using 
\beq\label{sharp}
\sharp \{Q\in \Bbb N^n:|Q|=m\} =\binom{n+m-1}{h} \leq 2^{n+m-1}
\eeq
 we obtain for $|x|\leq \delta/4$ and $\delta \leq 2/A$:
$$|\tilde{R}_N(z,x)|\leq C2^{n-1}\sum_{j=0}^{N-1}(A|z|)^j\Gamma(1+js)
(2|x|/\delta)^{N-j}(1-2|x|/\delta)^{-1}\leq $$
$$\leq C2^{n+N} \delta^{-N}\Gamma (1+(N-1)s) \sum_{j=0}^{N-1}|z|^j |x|^{N-j}.$$
From this and (\ref{tale}) the assertion follows.
\qed

Let $r\in \Bbb N^n$ a nonzero multiindex and let $\rho : \Bbb C^{n}\rightarrow \Bbb C^{n+1}$ be the map defined to be $\rho(x)=(x^r,x)$.
\begin{defis}\cite{Ram-Mart2}[p.6]
Let $f\in \Bbb C[[x]]$ be a formal power series in $\Bbb C^n$ (resp. smooth function in some domain). We shall say that $f$ is $k$-summable (resp. a Gevrey function) in the monomial $x^r$ if $f\in \rho^*\Bbb C\{x\}\{z\}_k$, that is $f$ is the pull-back of a formal series (resp. smooth function) $g(z,x)$ which is $k$-summable (resp. a Gevrey function) in the variable $z$ in some sector, uniformly in $x$ on a polydisk.
\end{defis}
Let $k>0$ and let ${d}$ denote a direction in the complex plane.
We define the {\bf Borel transform of order $k$} in the direction ${d}$ as
$$
{\cal B}_k f(t):=\frac{1}{2i\pi}\int_{\gamma_k}f(z)e^{(t/z)^k}d(z^{-k}).
$$
Here we assume that $f$ is holomorphic in a sector $S=\{z\in {\Bbb C^*}:
|z|<\rho,|d-\arg z|< \alpha \}$ where $\rho >0$ and $\alpha>\pi/(2k)$ and $\gamma_k$
is the path from 0 along the ray $\arg z = d-\alpha_1$ till $|z|=\rho_1$ then
 along the circle $|z|=\rho_1$ to the ray $\arg z = d+\alpha_1$ and then back to 
the origin along this ray. Here $0<\rho_1<\rho$ and $\pi/(2k)<\alpha_1 <\alpha$. 
If $\hat f = \sum_{n=1}^{\infty}f_n z^n$ is a formal power series,
 then the formal Borel transform of order $k$ is defined as the power series
$$
{\cal B}_k \hat f(t) = \sum_{n=1}^{\infty}\frac{f_n}{\Gamma(n/k)} t^{n-k},
$$
where $\Gamma(x)$ is the Gamma-function.
We define the {\bf Laplace transform of order $k$} in the direction $d$
 (inverse Borel transform) as
$$
{\cal L}_k \tilde f(z):=\int_{0}^{\infty:d}\tilde f(t)e^{-(t/z)^k}d(t^{k}).
$$
An equivalent definition for $k$-summability in a direction $d$ is:
\begin{defis}\cite{ramis-ksum, ramis-panora}
A formal power series $\hat f$ at the origin of the complex plane will be said
to be $k$-summable in the direction ${d}$ if its formal $k$-Borel transform
 defines a holomorphic function in a neighborhood of the origin which can be
 continued holomorphically in some small sector bisected by the direction 
${d}$ and is of exponential growth of order at most $k$ at infinity. In
 this case, the function ${\cal L}_k \circ {\cal B}_k \hat f$ is holomorphic
 in a {\bf large sector} bisected by ${\bf d}$ and of opening $>\pi/k$.
 Moreover, this is the unique function which admit $\hat f$ as asymptotic
 expansion in this sector.
\end{defis}
Another equivalent definition has been given by Tougeron as follows: Let $\eta$, $R$ be positive numbers and let $k>1/2$.
Let ${\cal S}(\eta,R)$ denote the sector
$$
{\cal S}_{d, \pi/k}(\eta,R):=\left\{z\in\Bbb C^*\;|\;|\arg z-d| < \frac{\pi}{2k}+\eta, 0<|z|<R \right\}.
$$
Let $\theta>0$.
Let us define the {\it sectorial neighborhood of the origin} of order $q$ to be
$$
{\cal S}_{q,d,\pi/k,\theta}(\eta,R):=\left\{z\in\Bbb C\;|\;|z|<\frac{R}{(q+1)^{\theta}} \right\}\cup {\cal S}_{d, \pi/k}(\eta,R)
$$
\begin{theos}[Tougeron's definition of summability]\cite{tougeron-ens}
A function $f$ is a $k$-sum in the direction $d$ if and only if it has a representation as a sum $\sum_{q=0}^{+\infty}f_q$ of functions $f_q$, each of which is holomorphic in the sectorial neighborhood ${\cal S}_{d,\pi/k,1/k}(\eta,R)$ and satisfies
$$
\|f_q\|_{{\cal S}_{q,d, \pi/k,1/k}(\eta,R)}:=\sup_{z\in {\cal S}_{q, d, \pi/k,1/k}(\eta,R)}|f_q(z)|\leq C\rho^q,
$$
for some constants $\eta, R, C, \rho$ independent of $q$.
\end{theos}
\section{Main results}

Our first main result shows that a series of functions defines a Gevrey function on a sector if the Borel transforms satisfy good estimates in some well chosen domains.
 \par Let $k$ be a positive integer, $\alpha$ and $\beta$ real numbers with $\alpha <\beta$. Define for $0\leq\epsilon <(\beta -\alpha)/2, \rho>0$:
\beq\label{newS}
S_{\epsilon}(\rho)=\{z\in \Bbb C^*\mid \alpha+\epsilon\leq \arg z \leq 
\beta -\epsilon, |z|\leq \rho\}.
\eeq 

\begin{theos}\label{main1} Let $k, \alpha$ and $\beta$ be as above with
$\beta -\alpha>\pi/k.$ Let $\rho >0, R>0.$
 Suppose $g(z,y)=\sum_{Q\in \Bbb N^n} g_Q(z)y^Q$ is a scalar-valued 
holomorphic function in
 $S_0(\rho)\times \overline{P}_n(0,R)$. 
Moreover, 
suppose $({\cal B}_k g_Q)(t)$ exists and is holomorphic in 
$\overline{\Delta}_m:=\overline{P}_1(0,c m^{-\gamma})$ and satisfies 
\beq\label{small in r_m}
|({\cal B}_k g_Q)(t)|\leq K^m \mbox{ in }\overline{\Delta}_m,\mbox{ if }  m=|Q|\geq 1,
\eeq
where $\gamma\geq 0$ and $c$ and $K$ are positive constants.
\par Then for all $\epsilon \in (0,(\beta -\alpha)/2)$ the function $g(z,y)$ is a Gevrey function of order $\gamma +\frac{1}{k}$ in $z$
 in $S_{\epsilon}(\rho')$ 
uniformly in $y$ in $P_n(0,R')$ for some $\rho'\in (0,\rho)$ and $R'\in (0,R)$ 
both depending on 
$\epsilon$. 
Moreover, if $0\leq \gamma <\frac{\beta -\alpha}{\pi }-\frac{1}{k}$ then
 $g(z,y)$ is a $\frac{k}{k\gamma +1}$-sum.
\end{theos}
\begin{defis}
We shall say that the linear part $s=\sum_{i=1}^n \lambda_i x_i\partial/\partial x_i$ is diophantine of the type $\gamma\geq 0$ if there exists $c>0$ such that, for all $Q\in \Bbb N_2^n$, for all $1\leq i\leq n$ then 
\beq\label{dio}
|(Q,\lambda)-\lambda_i|>\frac{c}{|Q|^{\gamma}}\mbox{ unless }
(Q,\lambda)-\lambda_i=0.
\eeq
\end{defis}
Our second main result gives the Gevrey property of a sectorial normalizing tranformation of a well prepared vector field.
\begin{theos}\label{main2}
If the linear part of (\ref{champs-2}) is diophantine of type $\gamma \geq 0$ and
the assumptions of theorem \ref{main-stolo-classif} are satisfied, then 
the {\bf sectorial normalizing biholomorphisms} defined by theorem \ref{main-stolo-classif} are {\bf Gevrey functions of order $(1+\gamma)/k$} in the resonance monomial $x^r$. Moreover, if $\gamma=0$ then the formal normalizing transformation is $k$-summable. 
\end{theos}
\begin{rems}[Important remark]
Using theorem \ref{main1}, we could show that if $0\leq\gamma<1$, then the formal normalizing transformation is $\frac{k}{\gamma +1}$-summable. Nevertheless, we should emphasize that there is no $\lambda\in\Bbb C^n$ which satisfies $(\ref{dio})$ with $0\neq \gamma <1$. This was pointed out by Yann Bugeaud who refers to \cite{schmidt}.  
However, for a fixed non zero $\lambda\in \Bbb C^n$, there are infinite sequences of multiindexes $\{Q_m\}$ such that 
$|(Q_m,\lambda)-\lambda_i|>\frac{c}{|Q_m|^{\gamma}}$ with $0\neq \gamma <1$ unless $(Q_m,\lambda)-\lambda_i=0$ . Hence, if it happens that, in our normalization process, the sole monomials that appears in the Taylor expansions of our objects belong to such a sequence, then we will obtain the claimed summability property.

It is a remarkable fact that we obtain the summability property even when some singularities accumulate at the origin in the Borel plane. This is due to the slow rate ($\gamma <1$) at which this accumulation occurs. The fact that there are no $\lambda$ which satisfies $(\ref{dio})$ with $\gamma<1$ has nothing to do with this phenomenon. It is just an arithmetic property. 
\end{rems}
In fact, we show that $\phi_i(y)=g_i(y^r,y)$ where $g_i(z,y)$, given by $(\ref{exp g})$,  is shown to be a Gevrey function in $z$ in some sector at the origin of $\Bbb C$, uniformly in $y$ in a polydisk centered at the origin in $\Bbb C^n$.

As far as we know, it is the first time that such an interplay between the rate of accumulation of small divisors at zero and the Gevrey character of the normalizing transformation is characterized. In other situations, such an interplay seems to be guessed (see for instance \cite{simo-averaging, lochak-simult,lombardi-nf}).

We can show that the formal $k$-Borel transform of $g_{i,Q}(z)$ has no
 singularity in the disc centered at the origin and of radius 
$r<\inf_{|P|\leq |Q|}\{|(P,\lambda)-\lambda_i|\neq 0\}$. This is the main 
reason for which the $g_i$'s are $k$-sums when there are no small divisors 
(i.e. there exists $c>0$ such that for all $P\in \Bbb N_2$, $|(P,\lambda)-
\lambda_i|>c$ if the number on the left hand side is not zero). In fact, the
 ${\cal B}_k g_{i,Q}$'s are holomorphic on the {\it same} disc centered at the
 origin and of radius $c/2$ (it is easy to show that the $g_{i,Q}$ have
 asymptotic expansions $\hat{g}_{i,Q}$ which are $1/k$-Gevrey power series). 
In particular, this is the case in dimenson $2$ \cite{Ram-Mart1, Ram-Mart2}.
 Our main result will quantify the interplay between the small divisor 
phenomenon and the Gevrey character of the normalizing transformation.

As in the proof of theorem \ref{main-stolo-classif} \cite{Stolo-classif}
[p.132-135], the main theorem reduces to the proof of the Gevrey character 
(in the variable $z$, uniformly in the variables $x$) of the linearizing 
transformation of the associated non-linear system $(\ref{irreg})$ with 
an irregular singularity at the origin. The proof of this fact relies on 
theorem \ref{main1} which proof is postponed to the end of the article.
\section{Proof of theorem \ref{main2}}
In the same way as the proof of theorem \ref{main-stolo-classif} reduces to 
 the proof of the sectorial linearization of the non-linear system 
(\ref{irreg}), the proof of theorem \ref{main2} reduces to the proof of the
 Gevrey character  of the sectorial linearization of (\ref{irreg}) (see section \ref{sect-fn}).
It is sufficient to consider the case $\beta=1$ in (\ref{irreg}).
 We will consider a little bit more generally:
\beq\label{Irreg}
z^{k+1}\frac{dx}{dz}=(\Lambda +z^kA)x+z^kf(z,x),
\eeq
where $\Lambda=\mbox{diag}\{\lambda_1,\ldots ,\lambda_n\}, A=\mbox{diag}\{\alpha_1,\ldots ,\alpha_n\}$ and
 \[f(z,x)=\sum_{j\in \Bbb N_2^n}f_j(z)x^j,\]
a convergent series for $|x|\leq\rho_1, |z|\leq \rho_2$ and the coefficients
 $f_j(z)$ are $\Bbb C^n$-valued holomorphic functions for $|z|\leq \rho_2$.
\par We introduce the following hypotheses:
\begin{itemize}
\item $(H_1)$ the eigenvalues $\lambda_i$ are not all 0 and all have a
 nonnegative imaginary part,
\item $(H_2)$ if there are real eigenvalues $\lambda_i$ then either these are
 all 0 and then $\Re \alpha_i>0$ or 
there exists $\lambda_{i_0}\in \Bbb R^*$ such that
$$\min_{i\neq i_0}{\mbox{Re }\left(\alpha_i -\frac{\lambda_i}{\lambda_{i_0}}
\alpha_{i_0}\right)}>0,
$$
where the minimum is taken over all indices $i\neq i_0$ such that $\lambda_i
\in \Bbb R.$
\item $(H_3)$ if there exists a resonance relation $(Q,\lambda)=\lambda_i$
for some $Q\in \Bbb N_2^n$ then $\alpha_i-(Q,\alpha)\not\in \Bbb N.$   
\item $(H_4)$ if $\lambda_i$ is not real then $f_i/x_i$ is  holomorphic in a neighborhood of the origin.
\item $(H_5)$ there exists $r=(r_1,\ldots ,r_n)\in \Bbb N_1^n$ such that if 
${\cal I}=\{1\leq i\leq n\;|\;r_i\neq 0\}$  then
\begin{itemize}
\item $\forall i\in {\cal I}$, $f_i=x_i\tilde f_i$ and $\tilde f_i$ is holomorphic in a neighborhood of the origin.
\item $\sum_{i\in {\cal I}}r_i\tilde f_i =0 $.
\end{itemize}
\end{itemize}

\par Let
\begin{equation}\label{secteur} 
S_j(\epsilon,\rho)=\{z\in \Bbb C^*||\arg z -\frac{\pi}{k}(j+\frac{1}{2})|\leq
\frac{\pi}{k} -\epsilon, |z|\leq \rho\}
\end{equation}
where $\rho>0, j=0,\ldots ,2k-1$ and $0<\epsilon <\pi /k$ fixed.

\par Then 
\begin{theos}\label{ODEGevrey}
Assume that hypotheses $(H_1),(H_2),
(H_3)$ and $(H_4)$ are satisfied and that $\Lambda$ is diophantine
 of type $\gamma\geq 0$. Then for even $j$ with $0\leq j \leq 2k-2$  
there exists a unique change of variables $x=y+g^j(z,y)$ with $g^j(z,0)=0,
 D_yg^j(z,0)=0$ which transforms (\ref{Irreg}) into
\beq\label{Irreg0}
z^{k+1}\frac{dy}{dz}=(\Lambda +z^kA)y,
\eeq
and where $g^j(z,y)$ is a Gevrey function of order $(1+\gamma)/k$ in $z$ in 
 $S_j(\epsilon,\rho)$
 uniformly for $y\in \overline{P}_n(0,R)$ for all $\epsilon$ in $(0,\pi /(2k))$ and positive numbers $\rho$ and $R$ depending on $\epsilon$. 
\par \par If $\gamma =0$ then $g^j(z,y)$ is $k$-sum of $\hat{g}(z,y)$ in the 
direction $\frac{\pi}{k}(j+\frac{1}{2})$ where $\hat{g}(z,y)$ is the asymptotic expansion in
 $z$ of $g(z,y)$.
\par If all $\lambda_h$ are real then these statements also hold for odd $j$
with $1\leq j\leq 2k-1$
 and direction $\pi /2$ replaced by $-\pi/2.$
\par If condition $(H_5)$ is satisfied then $x^r=y^r$.
\end{theos}

Except for the Gevrey property this is theorem 2.7.1 in
 \cite{Stolo-classif}. However, the Gevrey property follows from theorem 
\ref{main1} once condition (\ref{small in r_m}) has been verified and then
 theorem \ref{ODEGevrey} follows. For simplicity, we shall verify this 
condition first in the case $k=1$. We shall indicate at the end of the 
section the changes in the proof in case $k>1$. Moreover, it is sufficient to 
prove the theorem only for the case $j=0$ since the other cases may be 
obtained from this by a rotation of the independent variable $z$.
\par In the next two subsections we derive
 properties of and estimates for the Borel transform of the coefficients $g_Q(z)$ 
in the expansion (\ref{exp g})
on $S_0\times\overline{P}_n(0,R)$. Here $g_Q$ is holomorphic on $S_0:=S_0(\epsilon,\rho)$ (cf. (\ref{secteur})).

\subsection{Properties of the Borel transform of $g_Q$ in case $k=1$}
 The function $g$ satisfies
\begin{equation}\label{3}
z^2\frac{d}{dz}{g(z,y(z))}=(\Lambda +zA)g(z,y(z))+zf(z,y+g(z,y(z)))
\end{equation}
together with
$$
z^2\frac{d y(z)}{dz} = (\Lambda +zA)y(z).
$$
From the series expansion for $g(z,y)$ and Cauchy's inequality it follows that
 there exists a positive constant $M$ such that
\beq\label{estgQ}
|g_Q(z)|\leq MR^{-|Q|}, z\in S_0=S_0(\epsilon,\rho).
\eeq
Let $g_{{\bf e}_l}:={\bf e}_l$ for $l=1,\ldots ,n$, where  ${\bf e}_l:=(\delta_{i,l})_{1\leq i\leq n}$ with $\delta_{i,l}=0$ if $i\neq l$ and $1$ otherwise. 
Then 
\[ f(z,y+g(z,y))=\sum_{j\in \Bbb N_2^n}f_j(z)(\sum_{Q\in \Bbb N_1^n}g_Q(z)y^Q)^j\]
and so
\begin{equation}\label{4}
f(z,y+g(z,y))=\sum_{Q \in \Bbb N_2^n}
y^Q t_Q(z)
\end{equation}
where
\begin{equation}\label{5}
t_Q(z):=\sum_{2\leq |j|\leq |Q|}f_j(z)\Sigma '\prod_{l=1}^n \prod_{q=1}^{j_l}
(g_{i_{l,q}}(z))_l, z\in S_0.
\end{equation}
Here $(v)_l$ denotes the $l$th component of a vector $v\in \Bbb C^n$,
 and $\sum'$ denotes
 that the sum has to be taken over all $i_{l,q}\in \Bbb N_1^n$ such that 
$\sum_{l=1}^n\sum_{q=1}^{j_l} i_{l,q}=Q$. So
$1\leq |i_{l,q}|\leq |Q|-|j|+1\leq |Q|-1$ and $t_Q=0$ if $|Q|\leq 1$.
From this and (\ref{3}) it follows that
for $Q\in \Bbb N^n_2$:
\begin{equation}\label{7}
z^2g_Q'(z) +(\lambda_Q+z\alpha_Q)g_Q(z)=zt_Q(z),
\end{equation}
where $\lambda_Q:=(\lambda,Q)-\Lambda, \alpha_Q:=(\alpha,Q)-A,\alpha=(\alpha_1,
\ldots ,\alpha_n).$ 
\par Let $w_Q:=z^{|Q|}g_Q$. Then it follows from (\ref{7}) and (\ref{5}) that
\begin{equation}\label{9}
z^2w_Q'(z) +(\lambda_Q+z\beta_Q)w_Q(z)=zu_Q(z),
\end{equation}
where $\beta_Q:=\alpha_Q-|Q|Id$ and
\[u_Q(z):=z^{|Q|}t_Q(z)=\sum_{2\leq |j|\leq |Q|}f_j(z)
\Sigma '\prod_{l=1}^n \prod_{q=1}^{j_l}(w_{i_{l,q}}(z))_l\] 
since in $\sum'$ we have 
$\sum_{l=1}^n\sum_{q=1}^{j_l} i_{l,q}=Q$.
\par Let $G_Q:={\cal B}g_Q$,
 $W_Q={\cal B}w_Q$ and $U_Q={\cal B}u_Q$. These functions exist in
 $S':=\{t\in \Bbb C^*|2\epsilon 
\leq \arg t \leq \pi -2\epsilon\}$ and in a neighborhood of the origin
since $\hat{g}_Q(z)$ is a series of Gevrey order 1. Moreover, $W_Q(t)$ and 
$U_Q(t))$ are $O(t^{|Q|-1})$ on $S'$ since 
$w_Q(z)=O(z^{|Q|})$. Furthermore
\beq\label{difW}
G_Q(t)=\frac{d^m}{d t^m}W_Q(t) \mbox{ where } m=|Q|.
\eeq
Since $w_{{\bf e}_l}=zg_{{\bf e}_l}=z{\bf e}_l$ and ${\cal B}z=1$ we have $W_{{\bf e}_l}=
{\bf e}_l$ for $l=1,\ldots ,n$.
\par We apply the Borel transform to both sides of (\ref{9}).
Since
${\cal{B}}(z^2\frac{d}{dz}w_Q(z))(t)=tW_Q(t)$, we obtain from (\ref{9})
\begin{equation}\label{9a}
(t+\lambda_Q)W_Q+\beta_Q*W_Q=1*U_Q
\end{equation}
where 
\begin{equation}\label{10}
U_Q=\sum_{2\leq |j|\leq |Q|}(f_j(0)+({\cal{B}}f)_j*)
\Sigma '\prod_{l=1}^n* \prod_{q=1}^{j_l}*(W_{i_{l,q}})_l.
\end{equation}
Here, $\prod *$ denotes the product with respect to the convolution product
 (see also \cite{costin-duke, braaksma-transeries}).
Differentiating (\ref{9a}) we get
\beq\label{difW'}
(t+\lambda_Q)W_Q'+(1+\beta_Q)W_Q=U_Q.
\eeq


\subsection{Proof of condition (\ref{small in r_m}) in case $k=1$} 
\par Define for $m\in \Bbb N_1$
\beq\label{defr_m}
\rho_m:=\min \{|(P,\lambda)-\lambda_j|\mid P\in \Bbb N_1^n, |P|\leq m,j=1,\ldots,n,
 (P,\lambda)\neq\lambda_j\}.
\eeq
The diophantine condition (\ref{dio}) 
implies that there exists a positive constant
$c$ such that 
\beq\label{dio1}
\rho_m\geq cm^{-\gamma}.
\eeq
First we prove 
\begin{lemms}\label{lem1}
 $U_Q(t)$ and $W_Q(t)$ are holomorphic for $|t|<\rho_m$ 
and have a zero of order at least $m-1$ at 0, where $m=|Q|$.
\end{lemms}
{\Proof} For $m=1$ we have $W_{{\bf e}_l}={\bf e}_l, U_{{\bf e}_l}=0.$
 Next suppose 
that the property holds for some $m\geq 1$.
 Let $P\in \Bbb N^n$ with $|P|=m+1.$
Since $t^{a-1}*t^{b-1}=B(a,b)t^{a+b-1}$ we may deduce from  (\ref{10}) 
that $U_P(t)$ is holomorphic for $|t|< \rho_m$ and has a zero of order $m$ at 0.
If $\lambda_{P,l}\neq 0$ then  (\ref{difW'}) implies that $W_{P,l}$ is
 holomorphic for $|t|<\rho_{m+1}$ and has a zero of order $m$ at 0, whereas if
 $\lambda_{P,l}=0$
then using $(H_3)$ we may construct a formal power series solution of the $l$th component of
 (\ref{difW'}) with a zero of order $m$ at 0 and this series converges for 
$|t|<\rho_m$ since 0 is a regular singular point of (\ref{difW'}). 
 \qed

\par From hypothesis $(H_2)$ we deduce 
\begin{lemms}\label{lem2}
 There exist $m_0\in \Bbb N_1$ and $\delta_0>0$ such that for all $l\in
\{1,\ldots ,n\}$
\beq\label{m0}
\Re ((\alpha,Q)-\alpha_l)\geq\delta_0 \mbox{ if } (\lambda, Q)=\lambda_l, |Q|\geq
m_0.
\eeq
\end{lemms}
\Proof 
We may order the $\lambda_j$ such that $\lambda_j
\in \Bbb R$ if $j\leq q$ and $\Im \lambda_j>0$ if $j>q$. Let $Q'$ arise from
 $Q$ by replacing the last $n-q$ components by 0 and let $Q'':=Q -Q'$. From $(Q,\lambda)=\lambda_l$ it follows that $\Im (Q'',\lambda)=
\Im \lambda_l$. Since $\min\{\Im \lambda_j| j= q+1,\ldots ,n\}>0$ if $q<n$  we
 see that there exist positive constants $M_1$ and $M_2$ such that 
\[|Q''|\leq M_1,|(Q',\lambda)|=|\lambda_l- (Q'',\lambda)|\leq M_2.\] 
\par First suppose that there is an index $i\leq q$ such that $\lambda_i\neq 0$. 
Let $\delta_1$ be the minimum of the lefthand side in hypothesis $(H_2)$.
Let $P=Q'-Q_{i_0} {\bf e}_{i_0},p_0=\Re \alpha_{i_0}/\lambda_{i_0}$.
 From hypothesis 
$(H_2)$ it follows that
$\Re (P,\alpha)
\geq |P|\delta_1+(P,\lambda)p_0$
and therefore $\Re (Q',\alpha)\geq |P|\delta_1+(P,\lambda)p_0 +Q_{i_0}\Re \alpha_{i_0}=|P|\delta_1+(Q',\lambda)p_0$. 
Hence
\[\Re (Q,\alpha)\geq |P|\delta_1 +(Q',\lambda)p_0+\Re (Q'',\alpha)\geq
|P|\delta_1-M_2|p_0|-M_1|\alpha|.\]
 Let $M_3>(M_2|p_0|+(M_1+1)|\alpha|)/
\delta_1.$ If $|P|<M_3$ then since $Q_{i_0}\lambda_{i_0}=(Q',\lambda)
-(P,\lambda)$ we get $|Q_{i_0}|< (M_2+M_3|\lambda|)/|\lambda_{i_0}|
=:M_4$ and so $|Q|< M_3+M_4+M_1=:M_0$. Therefore if $|Q|\geq M_0$ then $|P|\geq 
M_3$ and $\Re ((Q,\alpha)-\alpha_l)\geq M_3 \delta_1-M_2|p_0|-(M_1+1)|\alpha|
=:\delta_0 >0.$ 
\par Next suppose that $\lambda_i=0$ if $i\leq q$. Then $\Re \alpha_i >0$
for $i\leq q$. Let $\delta_2 =\min_{i\leq q} \Re \alpha_i$. So $\delta_2>0$.
It follows that $\Re \{(Q',\alpha)-\alpha_l\}\rightarrow \infty$ if
$|Q'|\rightarrow \infty$. Since $|Q''|\leq M_1$ we also have 
 $\Re \{(Q,\alpha)-\alpha_l\}\rightarrow \infty$ if
$|Q|\rightarrow \infty$.   \qed

\par We will give estimates of $W_Q$ in a neighborhood of 0 in terms of
\beq\label{defrho}
R_m(t):=\frac{|t|^{m-1}}{(m-1)!}, m\in \Bbb N_1.
\eeq 
\begin{lemms}\label{lem3}
There exist positive constants $K_0$ and $c_0$ with $c_0<1$ such that
for all $Q\in \Bbb N_1^n$
\beq \label{estWQ}
|W_Q(t)|\leq K_0^m R_m(t) \mbox{ if } m=|Q|,|t|\leq c_0 \rho_m.
\eeq
\end{lemms}
\Proof We choose $c_0\in (0,1)$
 as follows.  Let
$p_{Q,l}=|Q|^{-1}((Q,\alpha)-\alpha_l), l=1,\ldots ,n.$ Then there exists a 
constant $M>0$ such that $|p_{Q,l}|\leq M$ for all $Q\in \Bbb N_1^n$ and all
 $l\in\{1,\ldots ,n\}$. Now choose $c_0$ so small that
 $\Re \frac{pt+1}{t+1}> 0$ for all $|t|\leq c_0$ and all $p\in \Bbb C$ such that $|p|\leq M$.
\par From lemma \ref{lem1} it follows that there exist positive constants
$\mu_Q$ such that
\beq\label{12}
|W_Q(t)|\leq \mu_Q R_{|Q|}(t) \mbox{ if }|t|\leq c_0 \rho_{|Q|} 
\eeq
for all $Q\in \Bbb N_1^n$. We choose these constants first for all $Q$ with 
$|Q|<m_0$, where $m_0$ is given in lemma \ref{lem2}. For $|Q|\geq m_0$ we
 determine suitable $\mu_Q$ by means of a recurrence relation.
Suppose $\mu_Q$ have been determined for $1\leq |Q|<m$ such 
that (\ref{12}) holds for these $Q$. Here we assume that $m\geq m_0$.  
\par We first estimate $|U_Q|$ for $|Q|=m.$ There 
exists a positive constant $K_1$ such that 
\[|f_j(0)|\leq K_1^{|j|}, |({\cal B} f_j)(t)|\leq K_1^{|j|}\]
for all $t\in \Bbb C$ with $|t|\leq \rho_1$ and $j\in \Bbb N$. From (\ref{10})
 it follows that for all $|t|\leq c_0 \rho_{m-1}$ we have
\begin{equation}\begin{split}
 |U_Q(t)|\leq \sum_{2\leq |j|\leq |Q|}K_1^{|j|}|(1+1*)
\Sigma ''\prod_{q=1}^{|j|}*\{\mu_{i_q}R_{|i_q|}(t)\}|\leq \\
(1+\rho_1)\sum_{2\leq |j|\leq |Q|}K_1^{|j|}R_{|Q|}(t)
\Sigma ''\prod_{q=1}^{|j|}\mu_{i_q},
\end{split}\end{equation}
where in $\sum ''$ we sum over $i_q\in \Bbb N^n_1$ with
$\sum_{q=1}^{|j|}i_q=Q$.            
Using (\ref{sharp}) we obtain for all $|t|\leq c_0 \rho_{m-1}$ and $|Q|=m\geq 2$
\begin{equation}\label{13}
|U_Q(t)|\leq \nu_{Q}R_{m}(t)\mbox{ where }
\nu_{Q}:=c_1\sum_{h=2}^{|Q|}(2K_1)^h\Sigma ''\prod_{q=1}^h\mu_{i_q},c_1:=
(1+\rho_1)2^{n-1}.
\end{equation}
We will use this estimate in an integral representation of $W_Q$ for $|Q|\geq 
m_0$ which follows
 from  (\ref{difW'}):
\beq\label{intW}
W_{Q,l}(t)=(t+\lambda_{Q,l})^{-1-\beta_{Q,l}}\int_0^t
(s+\lambda_{Q,l})^{\beta_{Q,l}}U_{Q,l}(s) ds, l=1,\ldots ,n.
\eeq
 From now on we will fix $l\in\{1,\ldots ,n\}$
and delete the index $l$.
\par First we suppose $\lambda_Q=0.$ From (\ref{intW}) and (\ref{13}) we obtain
\[ |W_Q(t)|\leq \nu_Q \int_0^t|(\frac{s}{t})^{\beta_Q+1}\frac{s^{m-2}}{(m-1)!}
ds|\leq\nu_Q R_m(t) (\Re \alpha_Q)^{-1} \]
if $|t|\leq c_0\rho_{m-1}$. With lemma
 \ref{lem2} we see that (\ref{12}) holds with $\mu_Q\geq \delta_0^{-1} \nu_Q$. 
\par Next suppose $\lambda_Q\neq 0.$ Let $t=\lambda_Q \tau$
 and $v_Q(\tau)=U_Q(t)/R_{|Q|} (t)$. So $|v_Q(\tau)|
\leq\nu_Q$ if $|\tau|\leq |\lambda_Q|^{-1} c_0 \rho_{m-1}$.
 We subsitute $s=\lambda_Q \sigma$ in (\ref{intW}) with
$|Q|=m$. Since $\beta_Q=m(p_Q-1)$
we obtain for all $|t|\leq c_0 \rho_{m-1}$
\beq
|W_Q(t)|\leq \frac{|t|^{m-1}}{(m-1)!} \int_0^\tau
|(\frac{\sigma(1+\sigma)^{p_Q-1}}{\tau(1+\tau)^{p_Q-1}})^{m-1}
(\frac{1+\sigma }{1+\tau})^{p_Q-1}\frac{v_Q(\sigma)}{1+\tau} d\sigma|.
\eeq
Let $h(\sigma ):=\log \{\sigma(1+\sigma)^{p_Q-1}\}$, so that the first factor
 in the integrand equals $\exp \{(m-1)(h(\sigma)-h(\tau))\}$. Then 
$h'(\sigma)=\frac{1+p_Q\sigma }{\sigma (1+\sigma)}$. 
Due to the choice of $c_0$ made above we have
$\Re \frac{1+p_Q \sigma}{1+\sigma} >0$ for all $Q \in {\Bbb N}^n_{2}$ 
if $|\sigma|\leq c_0$.
Hence $\Re\{\frac{d}{d\xi}h(\xi \tau)\}=\Re \frac{1+p_Q \xi \tau}{\xi(1+\xi \tau)}>0$ if $|\tau|\leq c_0, 0<\xi \leq 1$ and consequently $\Re(h(\sigma)-
h(\tau))<0$ if $\sigma \in [0,\tau), |\tau|\leq c_0$ for all
 $Q\in {\Bbb N}^n_2$. It follows that there is a positive constant $c_2$
 independent of $Q$ and $t$ such that
\beq
|W_Q(t)|\leq R_m(t)\int_0^\tau 
|(\frac{1+\sigma }{1+\tau})^{p_Q-1}\frac{d\sigma}{1+\tau}|\nu_Q\leq c_2 \nu_Q
R_m(t)
\eeq 
if $|Q|=m, |t|\leq c_0\min \{|\lambda_Q|,\rho_{m-1}\}$, so in particular for $|t|
\leq c_0 \rho_m.$ Enlarging $c_2$ such that $c_2\geq \delta_0^{-1}$ we define $\mu_Q=c_2\nu_Q$ implying
 estimate (\ref{12}) for both cases $\lambda_Q=0$ and $\lambda_Q\neq 0$.
\par  Thus with
(\ref{13}) we obtain the recurrence relation
\beq\label{recu}
 \mu_Q=c_1c_2 \sum_{h=2}^{|Q|} (2K_1)^h\Sigma ''\prod_{q=1}^h\mu_{i_q} 
\mbox{ if } |Q|\geq m_0 
\eeq
with $\Sigma''$ as before.
\par Define formally
$$F(x):=\sum_{Q\in \Bbb N_1^n} \mu_Q x^Q, x \in \Bbb C^n.$$
Let $c_1c_2\sum_{h=2}^\infty (2K_1F(x))^h=\sum_{Q\in \Bbb N_1^n}\tilde{\mu}_Q 
x^Q$ formally.
 From (\ref{recu}) we may deduce that $\tilde{\mu}_Q=\mu_Q$ if $|Q|\geq m_0$.
Hence $F-c_1c_2(2K_1F)^2(1-2K_1F)^{-1}=F_1$ where $F_1$ is a polynomial in $x$
 of degree $<m_0$ and $F_1(0)=0.$ This quadratic equation for $F$ has a unique 
holomorphic solution $F=F_1+O(F_1^2)$ as $F_1\rightarrow 0$. 
So the formal series $F(x)$
converges in a neighborhood of the origin and consequently there exists
a constant $K_0>0$ such that $\mu_Q \leq K_0^{|Q|}$ for all $Q\in \Bbb N_1^N.$
 This implies (\ref{estWQ}).   
\qed

\begin{coros}\label{GQsmall}
Let $Q\in \Bbb N_1^n$ 
and $m=|Q|$. There exists a positive constant $K$ such that ${\cal {B}}g_Q=G_Q$ satisfies
\beq\label{difGsmall}
|G_Q(t)|\leq K^m \mbox{ if }|t|\leq \rho_m /2.
\eeq
\end{coros}
\Proof Let $t$ be as above and $C$ the circle with radius $ \rho_m/2$ and center $t$. From Cauchy's theorem, (\ref{difW}) and (\ref{estWQ}) 
it follows that 
$$ |G_Q(t)|\leq \frac{m!}{2\pi }\int_C |\frac{W_Q(t+s)}
{s^{m+1}}ds|\leq m! K_0^m\frac{(\rho_m)^{m-1}}{(m-1)!}
(\rho_m/2)^{-m}, $$
 which implies (\ref{difGsmall}).
\qed

 Hence (\ref{small in r_m}) holds. So  the assumptions of theorem \ref{main1} 
with $\alpha=-\pi/2$ and $\beta =3\pi/2$ are satisfied and
 theorem \ref{ODEGevrey} follows in case $k=1.$

\subsection {Proof in case $k>1$}
As before we only need to prove the Gevrey property of the
 normalizing transformation $g(z,y)$. 
It is sufficient to consider
 the case $j=0$ only.
\par With  the coefficients $g_Q$ in (\ref{exp g}) -which are holomorphic in 
$S_0$- we now associate
$w_Q(z)=z^{k|Q|}g_Q(z)$ and $u_Q(z)=z^{k|Q|}t_Q(z)$.
Now instead of (\ref{9}) we have 
$$z^{k+1}w_Q'(z) +(\lambda_Q+z^k\beta_Q)w_Q(z)=z^ku_Q(z),\mbox{ where }
\beta_Q=\alpha_Q -k|Q|Id.$$
\par Let the operator $\sigma_k$ be defined by $(\sigma_k \phi)(z)=
\phi(z^{1/k})$ and define $\tilde{w}_Q=\sigma_k w_Q,\tilde{g}_Q=\sigma_k g_Q,
\tilde{u}_Q=\sigma_k u_Q$. Then 
\beq\label{tilde}
\tilde{w}_Q=z^{|Q|}\tilde{g}_Q,
k z^2 \tilde{w}_Q'(z) +(\lambda_Q+z\beta_Q)\tilde{w}_Q(z)=z\tilde{u}_Q(z).
\eeq
Let $\tilde{G}_Q={\cal B}\tilde{g}_{Q},\tilde{W}_{Q}={\cal B}\tilde{w}_{Q},
 \tilde{U}_{Q}={\cal B}\tilde{u}_{Q}.$ Then instead of (\ref{difW'}) we now
 have 
$$ (kt +\lambda_Q)\tilde{W}_Q'(t)+(k+\beta_Q)\tilde{W}_Q(t)=\tilde{U}_Q(t),$$
and instead of lemma \ref{lem1} we now have that
 $\tilde{U}_{Q}$ and $\tilde{W}_{Q}$ are holomorphic functions of 
$t^{1/k}$ for $|t|< \rho_m/k$ and both are $O(t^{m-1})$ as $t\rightarrow 0$. 
In lemma \ref{lem3} we have to replace $W_Q$ by $\tilde{W}_Q$ and
the condition on $|t|$ becomes $|t|\leq c_0 \rho_m/k.$
\par Corollary \ref{GQsmall} now becomes:
\\ {\em The $k$-Borel transform ${\cal B}_k g_Q(t)=G_Q(t)$ is holomorphic for }
{\em $|t|\leq c_1 \rho_m^{1/k}$ and $|G_Q(t)|\leq K^m$ on this set. Here $c_1$}
{\em and $K$ are positive constants and $m=|Q|.$} 
\par \Proof From the definition of the Borel transform it follows that
$\tilde{G}_Q=\sigma_k G_Q$ and therefore $G_Q(t)=\tilde{G}_Q(t^k).$
Let $\rho_m':=c_0 \rho_m/k.$ If $|t^k|<\rho_m'$ then 
$$ \tilde{W}_Q(t^k)=\frac{1}{2\pi i}\int_C \frac{\tilde{W}_Q(s^k)}{s-t}ds$$
 where $C$ is the positively oriented circle $|s|=(\rho_m')^{1/k}$. From 
(\ref{tilde}) 
we deduce that if $|t|<\rho_m'$ then
$$ \tilde{G}_Q(t)= \frac{d^m}{dt^m}\tilde{W}_Q(t)= \frac{1}{2\pi i}\int_C \tilde{W}_Q(s^k)\frac{d^m}{dt^m}
(s-t^{1/k})^{-1}ds.$$
There exists $K_1>0$ such that 
$$|\frac{d^m}{dt^m}(s-t^{1/k})^{-1}| \leq m! K_1^m |s|^{-km-1}\mbox{ if }|t|
\leq |s|^k/2 $$ for all $m\in \Bbb N$ and $s\neq 0.$ This follows easily by 
substituting 
$t=s^k u$. Using this estimate and the $k$-version of lemma \ref{lem3} in the 
preceding integral we obtain
$|\tilde{G}_Q(t)|\leq m (K_0K_1)^m /\rho_m'$
 if $|t|\leq \rho_m'/2$ and since $G_Q(t)=\tilde{G}_Q(t^k)$ this implies the 
corollary.
 \qed
\par Hence the assumptions of theorem \ref{main1} are satisfied with
$\alpha=-\pi/(2k)$ and $\beta=3\pi/(2k)$
if we replace 
$\gamma$ by $\gamma /k$
 and theorem \ref{ODEGevrey} follows in case $k>1.$
 
\subsection{Proof of theorem \ref{main2} from theorem \ref{ODEGevrey}}\label{sect-fn}
The argument is the same as in the proof of theorem 4.3.1 in \cite{Stolo-classif}. For the convenience of the reader we give the sketch of it. 
\par To a well prepared holomorphic vector field $X\in{\cal E}_{k,\lambda,\alpha}$ of the form $(\ref{champs-2})$ with $\beta=1$, we associated a germ of holomorphic vector field $\tilde X$ in $(\Bbb C^{n+1},0)$ 
$$
\tilde X= \sum_{i=1}^n{\left(x_i(\lambda_i+\alpha_{i}z^k)+z^kf_i(x)\right)\frac{\partial}{\partial x_i}}+z^{k+1}\frac{\partial}{\partial z}.
$$
It it tangent to the germ of variety $\Sigma=\{z=x^r\}$ at the origin and its restriction to it is equal to $X$. To $\tilde X$, we associate a non-linear system with irregular singularity at the origin
\begin{equation}\label{nonlin-sing}
z^{k+1}\frac{dx_i}{dz}=x_i(\lambda_i+\alpha_iz^k) + z^kf_i(x)\;\;\;\;i=1,\ldots,n.
\end{equation}
If the assumptions $(H_i)_{i=1,\ldots, 4}$ are satisfied then we can apply theorem \ref{ODEGevrey}.
Since the original vector field is well prepared, system $(\ref{nonlin-sing})$ satisfies condition $(H_5)$. Hence, the sectorial linearizing diffeomorphisms $x=y+g^j(z,y)$ preserve the monomial $x^r$. Therefore, we can restrict the associated vector field of $\Bbb C^{n+1}$ as well as the sectorial diffeomorphisms to $\Sigma$. We obtain that $x=y+g^j(y^r,y)$ transform $X=\tilde X_{|\Sigma}$ to 
$$
\sum_{i=1}^n{y_i\left(\lambda_i+\alpha_{i}(y^r)^k\right)\frac{\partial}{\partial y_i}}
$$
and have the good Gevrey properties. As already noticed in remark 4.3.1 of \cite{Stolo-classif}, we can associate different $\tilde X$ to $X$. Each of them differs from the others by a vector field vanishing on $\Sigma$ and gives rise to a set of linearizing sectorial diffeomorphisms. The point is that, althought these sets depends on the chosen $\tilde X$, their restriction to $\Sigma$ don't. Hence, the Gevrey property makes sense.

\section{Proof of theorem \ref{main1}}
\subsection{Estimates for $G_Q^{(N)}$}
Let $\rho_m:=c m^{-\gamma}$ for $m\in \Bbb N_1$ and $G_Q(t)={\cal B}_kg_Q(t)$. From the properties of 
$g(z,y)$ we deduce that (\ref{estgQ}) holds on $S_0:=S_0(\rho)$ 
(cf. (\ref{newS})) and that $G_Q$ is
 holomorphic on $S_1:=\{t\in \Bbb C^*|\alpha +\frac{\pi}{2k}+\epsilon_1 
\leq \arg t \leq \beta -\frac{\pi}{2k} -\epsilon_1\}$. Here we choose
$\epsilon_1$ sufficiently small: $0<\epsilon_1<\pi/(4k), \epsilon_1<
(\beta -\alpha-\pi/k)/2.$
 Since $G_Q(t)$ is holomorphic for $|t|\leq \rho_m$ and satisfies 
(\ref{small in r_m}), Cauchy's inequality shows that
\beq\label{GQsmal}
|\frac{G_Q^{(N)}(t)}{N!}|\leq K^m (\rho_m/2)^{-N}\mbox{ if }|t|\leq \rho_m/2.
\eeq
Next we estimate $G_Q(t)$ for $t\in S_1$. For these values of $t$ we have  
$$G_Q(t)=\frac{1}{2\pi i}\int_\Gamma e^{t^ks}g_Q(s^{-1/k})ds$$ 
where $\Gamma$ is a contour consisting of the arc $\Gamma_0$ 
of the circle 
$|s|=\rho^{-k}$ in the sector $|\arg (t^k s)|\leq \pi/2 +k\epsilon_1 $ and the
two half lines $\Gamma_\pm$ consisting of the part of the rays $\arg(t^k s)=
\pm (\frac{\pi}{2}+k\epsilon_1)$ outside this circle. Hence
\beq\label{intGQN} 
\frac{G_Q^{(N)}(t)}{N!}= \frac{1}{2\pi i}\int_\Gamma 
\frac{1}{N!}\{\frac{d^N}{d t^N}e^{t^ks}\}g_Q(s^{-1/k})ds. 
\eeq
Here
$$|\frac{1}{N!}\frac{d^N}{d t^N}e^{t^ks}|\leq\frac{1}{2\pi}
\int_{|\sigma|=\mu}
|\frac{\exp \{t^k(1+\sigma)^ks\}}{t^N\sigma^{N+1}} d\sigma|.$$
We choose $\mu\in (0,1)$ such that $|\arg (1+\sigma)|\leq \epsilon_1 /2$ for 
$|\sigma|=\mu$. Thus we obtain for $|\frac{(\mu t)^N}{N!}\frac{d^N}{d t^N}e^{t^ks}|
$ the upper bound $\exp (|t|(1+\mu)/\rho)^k$ for $s\in \Gamma_0$ and the upper
 bound $\exp \{-|t|^k(1-\mu)^k |s|\sin \frac{k\epsilon_1}{2}\}$ for $s\in 
\Gamma_\pm$ since on $\Gamma_\pm$ we have $\frac{\pi + k\epsilon_1}{2}\leq
|\arg (s(t(1+\sigma))^k)|\leq \pi$. 
 Using these estimates and (\ref{estgQ}) in (\ref{intGQN}) we deduce
 that for $t\in S_1$ we have
$$|\int_{\Gamma_0}|\leq M_1R^{-m}(\mu |t|)^{-N}
\exp (|t|(1+\mu)/\rho)^k\},$$
$$|\int_{\Gamma_\pm}|\leq M_2R^{-m}(\mu |t|)^{-N}\int_0^\infty
\exp \{-|t|^k(1-\mu)^k |s|\sin \frac{k\epsilon_1}{2}\}d|s|$$
and thus
 \beq
|\frac{G_Q^{(N)}(t)}{N!}|\leq M_3R^{-m}(\mu |t|)^{-N}\{|t|^{-k}+
\exp (|t|(1+\mu)/\rho)^k\}
\eeq
Here $t\in S_1$ and $M_1,M_2,M_3$ are some positive constants. We use this
 estimate for $|t|\geq \rho_m/2$ and (\ref{GQsmal}) for $|t|\leq \rho_m/2$. Then we
 get for $t\in (S_1\cup \overline{\Delta} (\rho_m/2))$
\beq\label{estGQN}
|\frac{G_Q^{(N)}(t)}{N!}|\leq M_4R_1^{m}R_2^N \rho_m^{-N}
\exp (|t|R_3)^k
\eeq
where $M_4, R_1, R_2$ and $R_3$ are positive constants.


\subsection{Estimates for $g_Q(z)$ in a subsector of $S$}
We use the Laplace representation for $g_Q$: 
\[g_Q(z)=g_Q(0)+\int_0^{\infty :\theta}e^{-(t/z)^k}G_Q(t)dt^k=
g_Q(0)+ z^k \int_0^{\infty :k(\theta -\arg z)}e^{-s}G_Q(zs^{1/k})ds\]
for 
$\arg z\in [\alpha +2\epsilon_1, \beta-2\epsilon_1]$ and $|z|$ sufficiently 
small where $\alpha+\frac{\pi}{2k} +\epsilon_1\leq 
\theta\leq \beta -\frac{\pi}
{2k} -\epsilon_1, |\theta-\arg z|\leq \frac{\pi}{2k} -\epsilon_1.$ From this we
 deduce
\beq\label{intgN}
\frac{1}{N!}g_Q^{(N)}(z)=\int_0^{\infty :k(\theta -\arg z)}e^{-s}G(z,s)ds
\eeq
where
$$G(z,s)=\sum_{l=0}^k {k\choose l}z^{k-l}\frac{s^{(N-l)/k}}{(N-l)!}
G_Q^{(N-l)}(zs^{1/k}).$$

With (\ref{estGQN}) we obtain
$$\frac{1}{N!}|g_Q^{(N)}(z)|\leq M_4R_1^m\sum_{l=0}^k {k\choose l}|z^{k-l}|
R_2^{N-l}\rho_m^{l-N}
\int_0^{\infty :k(\theta -\arg z)}|e^{-s+(|z|R_3)^k|s|}s^{(N-l)/k}ds|.$$
As $|\theta -\arg z|\leq \frac{\pi}{2k}-\epsilon_1$ we have $\Re s\geq |s|\sin
(k\epsilon_1)$ on the path of integration.
 Let $|z|\leq (\frac{\sin (k\epsilon_1)}{2})^{1/k}
/R_3=:\rho'$. Then $|(zR_3)^ks|-\Re s\leq
 -(|s|\sin(k\epsilon_1))/2$ and therefore
\[\frac{1}{N!}|g_Q^{(N)}(z)|\leq
M_4 R_1^m \sum_{l=0}^k {k\choose l}|z^{k-l}| R_2^{N-l}
 \rho_m^{l-N}\int_0^\infty e^{-(|s|\sin (k\epsilon_1))/2}
|s|^{(N-l)/k})d|s|.\]
From this we deduce
\beq\label{estgN}
\frac{1}{N!}|g_Q^{(N)}(z)|\leq C_1^m (C_2/\rho_m)^N \Gamma (\frac{N}{k}+1)
\eeq
for  $z\in S_{2\epsilon_1}(\rho')$ (cf. (\ref{newS})).
Here $C_1$ and $C_2$
 are positive constants.

\subsection{End of proof of theorem \ref{main1}}
Since $\frac{\partial^N }{\partial z^N}g(z,y)=\sum_{Q\in \Bbb N_1^n} g_Q^{(N)}(z)y^Q$ we obtain
with the help of (\ref{estgN})
\beq\label{gzy}
 \frac{1}{N!}|\frac{\partial^N }{\partial z^N}g(z,y)|\leq\sum_{Q\in \Bbb N_1^n}
\Gamma (\frac{N}{k}+1) C_1^{m}(C_2/ \rho_m)^{N}|y|^m,
\eeq
provided the righthand side converges and where $m=|Q|$ and 
$z\in S_{2\epsilon_1}(\rho')$. Now
 use $(\ref{sharp})$ and 
the diophantine condition (\ref{dio1}) in (\ref{gzy}) and obtain for $z\in 
S_{2\epsilon_1}(\rho')$
\beq\label{g(z,y)N}
|\frac{1}{N!}\frac{\partial^N }{\partial z^N}g(z,y)|\leq \Gamma (\frac{N}{k}+1)
 C_2^{N}c^{-N}2^{n-1}\sum_{m=1}^\infty  m^{\gamma N}
(2C_1|y|)^m.
\eeq
Finally we use that for $0<\delta <1$ there exist positive constants 
$C_3$ and $C_4$ 
such that for $ |x|\leq 1-\delta$ and
 $\mu>0$ 
\beq\label{sum}
|\sum_{m=1}^\infty m^\mu x^m| \leq C_3C_4^\mu \Gamma(\mu)
\eeq
(cf. proof below). So the series in the righthand side of (\ref{gzy}) converges
for $|y|\leq R':=\frac{1-\delta}{2C_1}.$ Combining (\ref{g(z,y)N})  with
 (\ref{sum}) and Stirling's formula we see that there exists a positive
 constant $C_5$ such that
\beq\label{final}
|\frac{1}{N!}\frac{\partial^N }{\partial z^N}g(z,y)|\leq C_5^N\Gamma((\gamma
+\frac{1}{k})N)
\eeq 
if $|y|\leq R', z\in S_{2\epsilon_1}(\rho')$. This proves the Gevrey property in theorem \ref{main1} for sufficiently small positive
$\epsilon$ and this suffices.
\par Proof of (\ref{sum}): Let the lefthand side of (\ref{sum}) be denoted by
$F(\mu,x)$ where $|x|< 1$. We apply Hankel's formula
\[ \frac{m^\mu}{\Gamma (\mu +1)}=\frac{1}{2\pi i}\int_\Gamma e^{ms}s^{-\mu-1}ds
\]
to get
\[ F(\mu,x)=|\frac{\Gamma(\mu+1)}{2\pi i}\int_\Gamma s^{-\mu -1} 
\sum_{m=1}^\infty (e^s x)^mds|=|\frac{\Gamma(\mu+1)}{2\pi i}\int_\Gamma 
s^{-\mu -1}\frac{e^s x}{1-e^s x} ds|.\]
Here $\Gamma$ is a contour from $\infty e^{-\pi i}$ to $\infty e^{\pi i}$ 
turning once around 0 in the positive sense such that $|e^s x|<1$ on $\Gamma$. If
 $|x|\leq 1-\delta$  we may choose $\Gamma$ such that it consists of the
 circle $|s|=\frac{-\log (1-\delta}{2}$ and the parts of the half lines
 $\arg s=\pm\pi$ outside this circle. Then it is easy to estimate the 
integral by $2\pi C_3 C_4^\mu /\mu$ with some positive constants $C_3$ and
 $C_4$ and (\ref{sum}) follows.

\bibliographystyle{alpha}

\end{document}